\documentclass[12pt]{article}
\usepackage{geometry} 
\title{Subtraction Division Games}
\author{Elizabeth J. Kupin\thanks{Rutgers University Mathematics Department. \texttt{ekupin@math.rutgers.edu}}}
\date{} 

\usepackage{amsmath}
\usepackage{amssymb}
\usepackage{amsthm}
\usepackage{tikz}
\usetikzlibrary{matrix,arrows}
\usepackage{url}

\newtheorem{theorem}{Theorem}
\newtheorem{lemma}{Lemma}

\begin{document}

\maketitle

\section{Introduction}

A Subtraction-Division game is a two player combinatorial game with three parameters: a set $S$, a set $D$, and a number $n$. The game starts at $n$, and is a race to say the number 1. Each player, on their turn, can either move the total to $n-s$ for some $s \in S$ or to $\lceil\frac{n}{d}\rceil$, for some $d \in D$. This was first introduced as a puzzle in the MSRI Gazette \cite{first}, for the special case where $S=\{1\}$ and $D=\{2\}$. The case when $S=\{1,2,\ldots t-1\}$ and $D=\{t\}$ has been studied by Aviezri Fraenkel in \cite{Fraenkel1} and \cite{Fraenkel2}, as well as Alan Guo in \cite{Guo}, under the name $t$-Mark. We will focus on the cases where $|S|=|D|=1$, and investigate  which sets of parameters permit a winning strategy for the first player. 

In general, we will let $S=\{a\}$ and $D=\{b\}$, and think of these as being fixed while $n$ varies. We are then interested in classes of values $(a,b)$, for which we can characterize what $n$ will be wins or losses. To do this, we will consider the Sprague-Grundy value of the game. The Sprague-Grundy function, which we will denote $SG_{a,b}(n)$, or just $SG(n)$ if the values of $a$ and $b$ are clear from context, is 0 if the position is a first player loss, and non-zero if it is a first player win. 

$SG_{a,b}(n)$ is defined recursively, based on the moves available to the first player. To make this formal, we consider a digraph where each node is a state in the game (a particular value of $n$, for this game), and directed edges indicate the moves available from each state. We start by setting the winning state: $SG_{a,b}(1)=0$, since if the total starts at 1 and first player must move, there is no possible way for first player to get to 1 at the end of his turn. For all other states, we define the Sprague Grundy value to be the minimum excluded value (mex) of the Sprague Grundy values of its out-neighbors. That is, for higher values of $n$, $SG_{a,b}(n)=\text{mex}\{ SG_{a,b}(n-a), SG_{a,b}(\lceil\frac{n}{b}\rceil) \}$. A more general discussion of Sprague Grundy values is available here, and here. This also includes a proof of the claim that a non-zero value corresponds to a win, and a zero corresponds to a loss. We will develop only the properties of Sprague Grundy values that we will use.

Our game has the property that any state has out-degree at most $|S|+|D|$, and for most of the games we consider $|S|+|D|=2$. From this it follows that $SG_{a,b}(n)$ is bounded, for all $a$, $b$. Furthermore, we see that no two adjacent states can have the same SG value.

\section{Characterization of the Game Sequences}

\begin{theorem}If $a=1$ and $b$ even, there is a complete characterization for when $SG_{a,b}(n)$ is zero.\end{theorem}

Suppose that the largest divisor of $a$ that is relatively prime to $b$ is 1, and $b$ is even. Then we have the following more general result:
\begin{theorem}If $a$ and $b$ as above, there is a complete characterization of when $SG_{a,b}(n)$ is zero.\end{theorem}

We will prove Theorem 1 by building up the following structural Lemmas:

\subsection{Proof Theorem 1}

\begin{lemma} Let $D=\{2d\}$. If $n$ in even, the game $\{1, 2d,n \}$ has a first player winning strategy. If $n \equiv 2d+1, 2d+3, \ldots 4d-1 \mod 4d$, then the game $\{1, 2d, n\}$ has no first player winning strategy.  \end{lemma}

\begin{proof}
We will first show that if $n$ is even, $SG_{1,2d}(n)$ has a first player winning strategy, by induction on $n$. 

Base Case: $n=2$ is clearly a winning game for the first player: he simply subtracts 1 from the total, and wins.

Assuming that $SG_{1,2d}(2k-2d)$ is non-zero, we first consider the game for any even $n$. For all of these games, one option is to move the total to $k$. If must be that either $\{1, 2d, k\}$ is a losing position (in which case, it follows trivially that all these games are wins), or that $\{1, 2d, k\}$ is also a winning position. If so, the game $\{1, 2d, 2dk-2d+1\}$ mist be a losing position for first player, as the only options are to move to $n=2dk-2d$ and $n=k$, both of which are winning positions, and so first player must hand second player a winning position. But if $n=2dk-2d+1$ is a losing position, then $n=2dk-2d+2$ must be a win, $n=2dk-2d+3$ must be a loss, and it will alternate until we reach $n=2dk$.

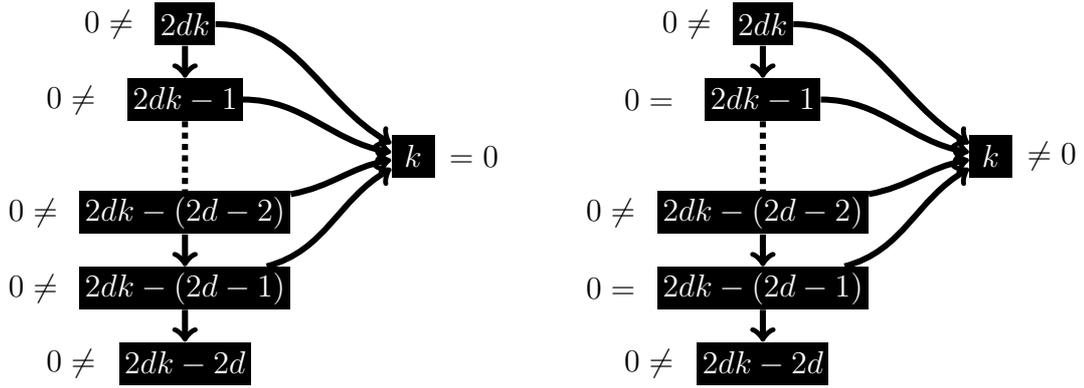
\begin{center}
\begin{tikzpicture}
\tikzstyle{vert}=[rectangle,fill=black, text=white, minimum size =16, inner sep=0]
\tikzstyle{label}=[rectangle,fill=white, text=black, minimum size =16, inner sep=0]
\draw (0,0) node[vert] (4dk){\,$2dk$\,};
\draw (0,-1) node[vert] (4dk-1) {\,$2dk-1$\,};
\draw (0,-2.5) node[vert] (4dk-2d+2){\,$2dk-(2d-2)$\, };
\draw (0,-3.5) node[vert] (4dk-2d+1){\,$2dk-(2d-1)$\, };
\draw (0,-4.5) node[vert] (4dk-2d){\,$2dk-2d$\,};
\draw (3,-1.75) node[vert] (2k){\,$k$\,};

\draw[->, line width= 0.08cm] (4dk) to (4dk-1);
\draw[line width=0.08cm, dotted] (4dk-1) --++ (4dk-2d+2);
\draw[->, line width= 0.08cm] (4dk-2d+2) to (4dk-2d+1);
\draw[->, line width= 0.08cm] (4dk-2d+1) to (4dk-2d);

\draw[->, line width= 0.08cm] (4dk) to[out=0, in=150] (2k);
\draw[->, line width= 0.08cm] (4dk-1) to[out=0, in=170] (2k);
\draw[->, line width= 0.08cm] (4dk-2d+2) to[out=10, in=190] (2k);
\draw[->, line width= 0.08cm] (4dk-2d+1) to[out=15, in=210] (2k);

\draw(-1,0) node[label](l1){$0 \neq $};
\draw(-1.5,-1) node[label](l1){$0\neq$};
\draw(-2,-2.5) node[label](l1){$0\neq $};
\draw(-2,-3.5) node[label](l1){$0\neq $};
\draw(-1.5,-4.5) node[label](l1){$0\neq$};
\draw(3.8,-1.75) node[label](l1){$= 0$};


\tikzstyle{vert}=[rectangle,fill=black, text=white, minimum size =16, inner sep=0, xshift=7.6cm]
\tikzstyle{label}=[rectangle,fill=white, text=black, minimum size =16, inner sep=0, xshift=7.6cm]
\draw (0,0) node[vert] (4dk){\,$2dk$\,};
\draw (0,-1) node[vert] (4dk-1) {\,$2dk-1$\,};
\draw (0,-2.5) node[vert] (4dk-2d+2){\,$2dk-(2d-2)$\, };
\draw (0,-3.5) node[vert] (4dk-2d+1){\,$2dk-(2d-1)$\, };
\draw (0,-4.5) node[vert] (4dk-2d){\,$2dk-2d$\,};

\draw (3,-1.75) node[vert] (2k){\,$k$\,};

\draw[->, line width= 0.08cm] (4dk) to (4dk-1);
\draw[line width=0.08cm, dotted] (4dk-1) --++ (4dk-2d+2);
\draw[->, line width= 0.08cm] (4dk-2d+2) to (4dk-2d+1);
\draw[->, line width= 0.08cm] (4dk-2d+1) to (4dk-2d);
\draw[->, line width= 0.08cm] (4dk) to[out=0, in=150] (2k);
\draw[->, line width= 0.08cm] (4dk-1) to[out=0, in=170](2k);
\draw[->, line width= 0.08cm] (4dk-2d+2) to[out=10, in=190] (2k);
\draw[->, line width= 0.08cm] (4dk-2d+1) to[out=15, in=210] (2k);

\draw(-1,0) node[label](l1){$0 \neq $};
\draw(-1.5,-1) node[label](l1){$0=$};
\draw(-2,-2.5) node[label](l1){$0\neq $};
\draw(-2,-3.5) node[label](l1){$0= $};
\draw(-1.5,-4.5) node[label](l1){$0\neq$};

\draw(3.8,-1.75) node[label](l1){$\neq 0$};

\draw(0,-6) node[rectangle,fill=white, text=black, text width=6cm](caption){If $SG(k)=0$, it's easy to see that all of $SG(2kd)$ through $SG(2kd-2k+1)$ will be zero.};
\draw(7.6,-6) node[rectangle,fill=white, text=black, text width=6cm](caption){If $SG(k)$ non-zero, then we will get alternation of the $SG$ values between zero and non-zero.};

\end{tikzpicture}
\end{center}

Given that $\{1, 2d, 2k\}$ is always a win, it's easy to see that if $n=4dk+2d+1, 4dk+2d+3, \ldots 4dk+4d-1$, first player's only moves in the game $\{1, 2d, n\}$ are to an even total and so to winning positions. This guarantees that second player will be handed a winning position, making this a losing game for first player.

\end{proof}

\begin{lemma} If $n\equiv 1, 3, 5, \ldots 2d-1 \mod 4d$, the base $2d$ representation of $n$ indicates whether or not the game $\{1, 2d, n\}$ is a first player win or loss. Consider the block of even digits immediately preceding the final (least significant) odd digit of $n$. If that block is even then $n$ is a loss,  and if it is odd then $n$ is a win. \end{lemma}
\begin{proof}
First we notice that, given Lemma 1 above, subtracting 1 from the total will always hand second player a winning position. First player's only hope for a win is to divide by $2d$, and move the total from $n=4dk+1$ to $n=\lceil\frac{4dk+1}{2d}\rceil = 2k+1$. If $n=2k+1$ is a loss, then $n=4dk+1$ will be a win, and vice versa. In fact, we will keep going with this process until we get to a point where we have a number that is equivalent to $2d+1, 2d+3, \ldots 4d-1 \mod 4d$, and we are guaranteed that this is a loss (by Lemma 1). The condition that each step has the opposite value as the step before it guarantees that our sequence of wins or losses will alternate along this path. Therefore, since we end in a loss, if we have taken an even number of steps to get to this point we have a loss, but if we have taken an odd number of steps we will have a win.

\bigskip

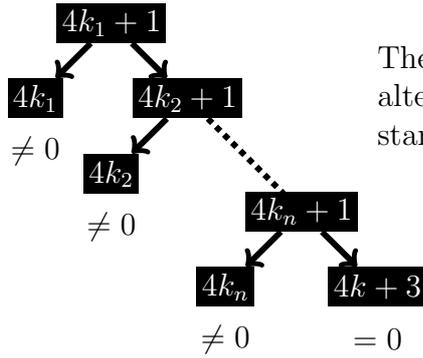
\begin{center}
\begin{tikzpicture}
\tikzstyle{vert}=[rectangle,fill=black, text=white, minimum size =15, inner sep=0]
\draw (0,0) node[vert] (4k+1){\,$4k_1+1$\,};
\draw (-1,-1) node[vert] (4k) {\,$4k_1$\,};
\draw (1,-1) node[vert] (2k+1){\,$4k_2+1$\,};
\draw (3.5,-3.5) node[vert] (end){\,$4k+3$\,};

\draw (-1,-1.7) node[rectangle, fill=white] (){$\neq 0$};

\draw (0,-2) node[vert] (2k){\,$4k_2$\,};
\draw (2.5,-2.5) node[vert] (k+1){\,$4k_n+1$\,};
\draw (1.5,-3.5) node[vert] (4kn){\,$4k_n$\,};
\draw (1.5,-4.2) node[rectangle, fill=white] (){$\neq 0$};
\draw (3.5,-4.2) node[rectangle, fill=white] (){$= 0$};

\draw (0,-2.7) node[rectangle, fill=white] (){$\neq 0$};

\draw (7,-1) node[rectangle, fill=white, text=black, text width=7cm](){The $SG$ values along the spine will alternate between zero and non-zero, starting with non-zero at the bottom.};

\draw[->, line width=0.08cm] (4k+1) to (4k);
\draw[->, line width=0.08cm] (4k+1) to (2k+1);
\draw[->, line width=0.08cm] (2k+1) to (2k);
\draw[->,line width=0.08cm] (k+1) to (end);
\draw[->,line width=0.08cm] (k+1) to (4kn);
\draw[dotted, line width=0.08cm] (2k+1) to (2.3,-2.3);
\end{tikzpicture}
\end{center}

The number of steps needed can be quickly computed from the  base $2d$ representation of $n$. Since $n\equiv 1, 3, \ldots 2d-1 \mod 4d$, the last (least significant) digit will be odd, and the next to last digit will be even. There is a block of even digits immediately preceding the final odd digit, with length at least one. Every time we divide by $2d$ and round up, we cut off the last digit, and add one to the (new) least significant digit. This has the net effect of shortening the length of this block of even digits by exactly one. Eventually, we reach a smaller number whose base-$2d$ representation ends in two odd numbers. This number must be equivalent to one of $2d+1, 2d+3, \ldots \text{ or } 4d-1 \mod 4d$, and so we know by Lemma 1 that this is a losing position for first player. The length of the block is exactly how many times we must divide by $2d$ to reach a number that is guaranteed to be a loss, and so the parity of that length tells us whether or not $n$ is a win.

\end{proof}

These two lemmas together prove Theorem 1 by giving a complete characterization.

When we start to look at subtracting numbers other than 1, we see a new type of pattern arise.

\subsection{Stuttering}

\begin{theorem} Let $a'$ be the largest divisor of $a$ that is relatively prime to $b$. Then the SG sequence for $a$ and $b$ will experience stuttering of length $s=\frac{a}{a'}$\end{theorem}

We will prove this by first showing the persistence of stuttering, once it occurs, and then by showing that stuttering occurring once is inevitable. It is as if we are looking at a proof by induction out of order, with the inductive step first. The reason for this presentation is that, unlike a standard proof by induction, the proof of the `base case' of the inevitability of stuttering is much more involved than the inductive step.

\begin{lemma} Stuttering of length $g=\gcd(a,b)$ will persist. \end{lemma}
\begin{proof} Consider a block of length $g$, that is a set of terms with indices $kg,kg-1, \ldots kg-g+1$, for some integer $k$. Because $g$ divides $b$, $\lceil \frac{kg}{b}\rceil=\lceil \frac{kg-1}{b}\rceil=\cdots= \lceil \frac{kg-g+1}{b}\rceil$, that is, these terms all depend on the same value via division. Moreover, because $g$ divides $a$, the set of terms with indices $kg-a, kg-a-1, \ldots kg-a-g+1$ is also a block of length $g$.

Suppose we have stuttering occur in a block, that is for some $m$, $SG(mg)=SG(mg-1)=\cdots=SG(mg-g+1)$. Consider the subsequent block, $SG(mg+a), SG(mg+a-1), \ldots SG(mg+a-g+1)$. Each term is the mex of a set of size 2. One element in the set is always identical across all terms, because they all depend on the same value under division. In this case, the value that they depend on via subtraction is identical as well, by assumption. Therefore it follows that $SG(mg+a) = SG(mg+a-1)= \cdots= SG(mg+a-g+1)$, and by induction this will occur for every block of the form $SG(mg+ka), SG(mg+ka-1), \ldots SG(mg+ka-g+1)$.
\end{proof}

Note that the term `subsequent' refers to the dependence under subtraction, and does not imply adjacent. In fact, we may have $a>>g$, in which case each sequence of subsequent blocks will be quite sparse. To get stuttering, we would need stuttering to appear in some block for each of the residue classes modulo $\frac{a}{g}$. In the following Lemma, we show that stuttering will occur in any sequence of subsequent blocks.

\begin{lemma} Stuttering of length $g=\gcd(a,b)$ is inevitable. \end{lemma}
\begin{proof} 
We will create a digraph that models subsequent blocks of length of $g$. Each vertex in this digraph will be a triple of SG values: $(x,y,z)$. $x$ and $y$ are two sample values in the block, that we hope will eventually be equal, whereas $z$ is the SG value that the next block depends on, via division.

For example, when $a=4$ and $b=2$, one triple could represent $$\left(SG_{4,2}(19), SG_{4,2}(20), SG_{4,2}(12)\right).$$ The subsequent block, $(23, 24)$, will depend on $(19,20)$ via subtraction, and on 12 via division. So just from looking at one triple, we have all the information we need to tell us the SG values in the next block.

We will add edges as follows: 

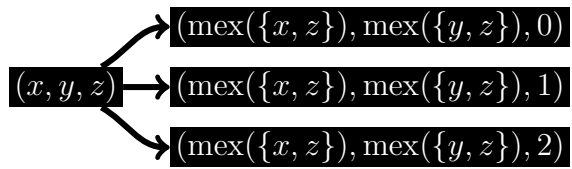
\begin{center}
\begin{tikzpicture}
\tikzstyle{vert}=[rectangle,fill=black, text=white, minimum size =15, inner sep=0]
\draw (0,0) node[vert] (base){\,$(x,y,z)$\,};
\draw (4,0.8) node[vert] (b0) {\,$(\text{mex}(\{x,z\}), \text{mex}(\{y,z\}), 0)$\,};
\draw (4,0) node[vert] (b1){\,$(\text{mex}(\{x,z\}), \text{mex}(\{y,z\}), 1)$\, };
\draw (4,-0.8) node[vert] (b2){\,$(\text{mex}(\{x,z\}), \text{mex}(\{y,z\}), 2)$\, };

\draw[->, line width= 0.08cm] (base) to [out=30,in=180](b0);
\draw[->, line width= 0.08cm] (base) to (b1);
\draw[->, line width= 0.08cm] (base) to [out=330, in=180](b2);
\end{tikzpicture}
\end{center}
Any sequence of blocks will correspond to a walk on this digraph. Fortunately, the digraph only has 27 vertices an 81 directed edges, so it is small enough to analyze by hand.

Here is a schematic of the interesting portions of the digraph. I have omitted vertices where $x=y$, because we know that once stuttering begins it will persist. Our primary concern, therefore, is how the graph behaves on the vertices where $x\neq y$, which reduces the number of vertices to 18. I have also depicted the following six vertices as sinks (and drawn them in gray): $(1,2,1)$, $(1,2,2)$, $(2,1,1)$, $(2,1,2)$, $(0,2,0)$, and $(2,0,0)$. For each of these, all three out edges point to vertices where $x=y$. For the first four $x=y=0$, and for the last two $x=y=1$. This further reduces the number of directed edges to consider, from 54 to 36.\\

\begin{tikzpicture}
\tikzstyle{vert}=[rectangle,fill=black, text=white, minimum size =16, inner sep=0]
\tikzstyle{sink}=[rectangle,fill=gray, text=white, minimum size =16, inner sep=0]
\draw (-1,0) node[vert] (201){\,$(2, 0,1)$\,};
\draw (3,0) node[vert] (021){\,$(0,2, 1)$\,};
\draw (-2,-1.5) node[vert] (202){\,$(2, 0, 2)$\,};
\draw (4,-1.5) node[vert] (022){\,$(0, 2, 2)$\,};
\draw (0,-2.25) node[sink] (020){\,$(0, 2, 0)$\,};
\draw (2,-2.25) node[sink] (200){\,$(2, 0, 0)$\,};
\draw (-2,-3) node[vert] (011){\,$(0, 1, 1)$\,};
\draw (4,-3) node[vert] (101){\,$(1, 0, 1)$\,};
\draw (-1,-4.5) node[vert] (012){\,$(0, 1, 2)$\,};
\draw (3,-4.5) node[vert] (102){\,$(1, 0, 2)$\,};

\draw (-1,-6.5) node[vert] (010){\,$(0,1,0)$\,};
\draw (3,-6.5) node[vert] (100){\,$(1,0, 0)$\,};
\draw (-0.25,-8) node[vert] (120){\,$(1, 2,0)$\,};
\draw (2.25,-8) node[vert] (210){\,$(2,1,0)$\,};

\draw (-2,-9.5) node[sink] (121){\,$(1,2,1)$\,};
\draw (0,-9.5) node[sink] (122){\,$(1,2,2)$\,};
\draw (2,-9.5) node[sink] (211){\,$(2,1,1)$\,};
\draw (4,-9.5) node[sink] (212){\,$(2,1,2)$\,};

\draw[->, line width= 0.08cm, color=gray] (201) to (020);
\draw[->, line width= 0.08cm] (201) to [out=10, in=170](021);
\draw[->, line width= 0.08cm] (201) to (022);

\draw[->, line width= 0.08cm, color=gray] (021) to (200);
\draw[->, line width= 0.08cm] (021) to [out=190, in=350] (201);
\draw[->, line width= 0.08cm] (021) to (202);

\draw[->, line width= 0.08cm] (202) to [out=210, in=160](010);
\draw[->, line width= 0.08cm] (202) to [out=285, in=75] (011);
\draw[->, line width= 0.08cm] (202) to [out=210, in=180](012);

\draw[->, line width= 0.08cm] (022) to [out=330, in=20] (100);
\draw[->, line width= 0.08cm] (022) to [out=285, in=75] (101);
\draw[->, line width= 0.08cm] (022) to [out=330, in=0](102);

\draw[->, line width= 0.08cm, color=gray] (011) to [out=0, in= 210]  (200);
\draw[->, line width= 0.08cm] (011) to [out=150, in=180] (201);
\draw[->, line width= 0.08cm] (011) to [out=105, in=255] (202);

\draw[->, line width= 0.08cm, color=gray] (101) to [out=180, in= 330] (020);
\draw[->, line width= 0.08cm] (101) to [out= 30, in= 0](021);
\draw[->, line width= 0.08cm] (101) to [out=105, in=255] (022);

\draw[->, line width= 0.08cm] (012) to (100);
\draw[->, line width= 0.08cm] (012) to (101);
\draw[->, line width= 0.08cm] (012) to [out=10, in=170] (102);

\draw[->, line width= 0.08cm] (102) to (010);
\draw[->, line width= 0.08cm] (102) to (011);
\draw[->, line width= 0.08cm] (102) to [out=190, in=350] (012);

\draw[->, line width= 0.08cm] (010) to (120);
\draw[->, line width= 0.08cm, color=gray] (010) to (121);
\draw[->, line width= 0.08cm, color=gray] (010) to [out=240, in=140] (122);

\draw[->, line width= 0.08cm] (100) to (210);
\draw[->, line width= 0.08cm,color=gray] (100) to [out=300, in=40] (211);
\draw[->, line width= 0.08cm,color=gray] (100) to (212);

\draw[->, line width= 0.08cm] (210) to [out=190, in=350]  (120);
\draw[->, line width= 0.08cm,color=gray] (210) to (121);
\draw[->, line width= 0.08cm,color=gray] (210) to (122);

\draw[->, line width= 0.08cm] (120) to [out=10, in=170]  (210);
\draw[->, line width= 0.08cm,color=gray] (120) to (211);
\draw[->, line width= 0.08cm,color=gray] (120) to (212);

\draw[dashed](-5,-5.5) to (7,-5.5);

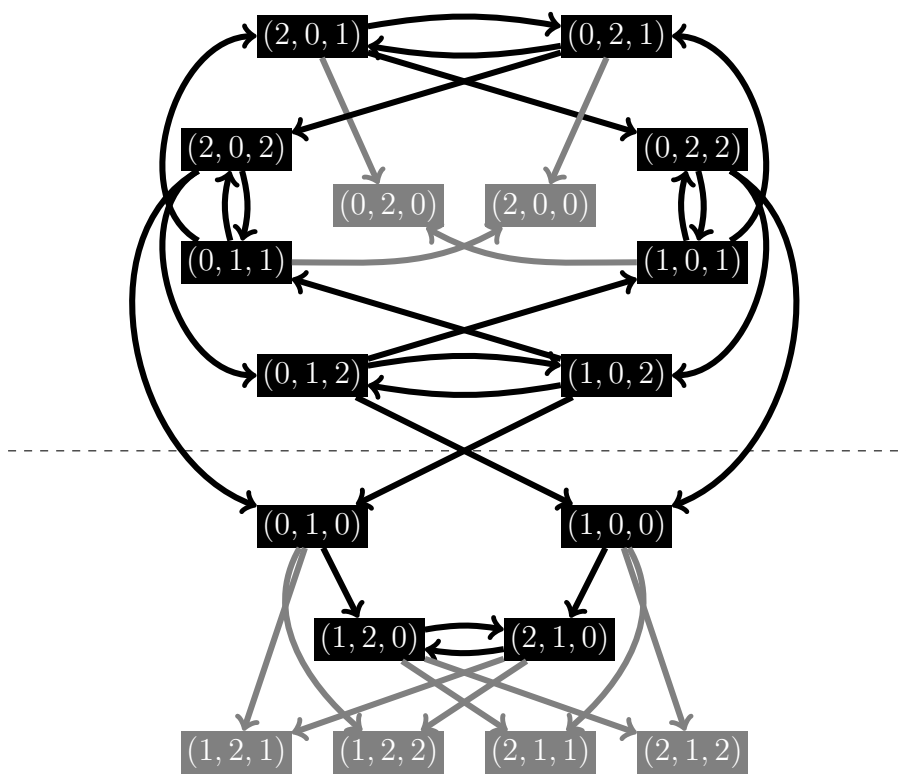
\draw(1,-11) node[rectangle, fill=white, text=black, text width =7.5 cm](){Digraph representing subsequent blocks.};

\end{tikzpicture}

We should be worried, since there are directed cycles where the first two terms are not zero, and an infinite walk around one of these cycles would correspond to an infinite sequence of blocks where stuttering does not appear. This cannot happen in practice, however.

Observe that above the dotted line all the vertices (other than the sinks) have $z\neq0$, and below the dotted line all the vertices (other than the sinks) have $z=0$. There are no arrows that go up from below the dotted line, so any directed cycle must stay either above or below the line. Therefore, an infinite walk without stuttering will have an infinite stretch with $z=0$, or an infinite stretch with $z \neq 0$. 

We will show in the following lemma that $z$ must vary between zero and non-zero values. This guarantees that wherever we begin, we will eventually pass to a sink or to a vertex below the dotted line. Moreover, once we are below the dotted line we must eventually move to a sink. Thus, stuttering is inevitable. 

\end{proof}

\begin{lemma}The value of $z$, above, cannot have an arbitrarily long sequence of zeros, or an arbitrarily long sequence of non-zero values. \end{lemma}
\begin{proof} 
When we consider the values of $z$, we see that this sequence depends on itself via subtraction: $z_i-a=z_j$, for some other value $j$. In fact, $j=i-d$. Therefore, it's impossible to have an infinite stretch of zeros. 

In theory, we could have an infinite stretch of non-zero values, but the only way this would be possible is if the value that $z_i$ depended on under division was zero for all $i$. That series will also depend on itself via subtraction, however, and so this less sparse sequence cannot be all zero. Therefore, the $z_i$ cannot be all zero or all non-zero.\end{proof}

We have only shown that stuttering of length $g=\gcd(a,b)$ will occur. To prove Theorem 3, we must show that we will get longer blocks. Armed with Lemmas 3 and 4, we can now prove Theorem 3:

\begin{proof}
Consider what happens in $SG_{a,b}(n)$, after stuttering has occurred. Far enough out in the sequence, the values in each block will be the same. It's natural, then, to think of the sequence not in terms of individual elements but in terms of the value of each block. 

The number of blocks `between' two subsequent blocks, that is how many blocks are skipped when we look from one block to the block it depends on via subtraction, has gone down by a factor of $g$. The block that one block depends on under division, however, hasn't changed. Therefore, far enough out in the sequence, the blocks have the same underlying digraph as the game \{$\frac{a}{g}$, $b$, $n$\}. 

The proof of the inevitability and persistence of stuttering tells us that if $\frac{a}{g}$ and $b$ have any common factors, we would expect to see stuttering of length $\gcd(\frac{a}{g},b)$. Since the elements are blocks, now, this would correspond to stuttering of length $g\cdot\gcd(\frac{a}{g},b)$ in the whole sequence. 

This process continues for as long as the amount we subtract by continues to have common factors with $b$. This will end once we reach $a^\prime$, the largest divisor of $a$ that is relatively prime to $b$. We will have stuttering of length $\frac{a}{a^\prime}$, and far out in the sequence, the underlying digraph that determines the values of the blocks will be the same as the digraph of the game $\{a^\prime, b, n\}$.
\end{proof}

\subsection{The Mis\'ere Game and Changed Initial Values}

The results from section 2.2 tell us that if all the prime factors of $a$ also appear as prime factors of $2d$ (i.e. $a^\prime=1$), then the game $\{a,b,n\}$ will eventually behave like the game $\{1,2d,n\}$. But it's not immediately clear how the characterization of the games $\{1, 2d, n\}$ can be applied to the game $\{a, 2d,n\}$. It turns out that we will need a more involved understanding of the structure of the simpler game, much of the details of which are postponed until Section 3. However in this section we will discuss how the two games are related, what results from above can be extended directly to these sequences, and how the upcoming results can be used to complete the characterization. 

Since we are studying a sequence that has eventual stuttering of length $a$, we will start by breaking the full sequence up into subsequences based on the residue classes mod $a$. That is, one sequence will be the terms with indices $a, 2a, 3a,$ etc., another will be the terms with indices $1, a+1, 2a+1, 3a+1,$ etc., and so on. Stuttering tells us that these sequences are equal for (subsequence) indices greater than $N$, so it doesn't matter which residue class we consider. For the rest of this section, we will look at only one (unspecified) residue class, and will take the view that our sequence of interest is made by taking $N-1$ fixed, arbitrary values, and then for all subsequent terms taking the recursive definition for the game $\{1,2d,n\}$: $$x_n= \text{mex} \left( x_{n-1}, x_{\lceil \frac{n}{2d} \rceil}\right).$$
 
We would like to prove an analog of Lemma 1 for this subsequence, that is, that even terms have a first player win and terms with indices in the range $n \equiv 2d+1, 2d+3, \ldots 4d-1 \mod 4d $ have no first player winning strategy. Of course, we have no control over the first $N-1$ terms, as they may be set arbitrarily (and in particularly, can perhaps be set in such a way that will influence subsequent terms), so the best we can hope for is that these results will hold for all indices above some threshold.

\begin{lemma} In a sequence of the game $\{1, 2d, n\}$, where the first $N-1$ terms have been set arbitrarily, there is a first player winning strategy for the even index terms starting at $2dN$. There is no first player winning strategy for terms with indices $n \equiv 2d+1, 2d+3, \ldots 4d-1 \mod 4d $, for all $n\geq 4d^2N-2d+1$. \end{lemma}

\begin{proof}
As in the proof of Lemma 1, we have an inductive argument for a winning strategy for even index terms. This essentially says that if $n=2d+2k$, then as long as $n=2d$ has a winning strategy, $n=2d+2k$ will also have a winning strategy. What we need to show is that the base case appears, no later than the term with index $2dN$. Suppose, then, that all even index terms from $N$ to $2dN-2$ have no winning strategy. The value at the term $2dN$ is determined by the mex over a set that includes the value at the term with index $N$, which is zero. Since the mex of a set that includes zero cannot be zero, the value at $2dN$ must be non-zero, corresponding to a winning strategy.

For the second part of the lemma, we depend on the first part. We know that for indices in that range, both elements in the mex will have even index. Once we are at a point where we are guaranteed that the indices are large enough, so that even index terms will be non-zero, we will have the mex over a set without zero, which must be zero. Since we will depend on numbers through both subtraction and division, this will not be guaranteed to occur until at least the index $4d^2N-2d+1$, as $\lceil \frac{4d^2N-2d+1}{2d}\rceil=2dN$. 
\end{proof}

Already, we see what type of characterization we may expect. If we have a term in the sequence given by $\{a, 2d, n\}$ and we wish to determine if there is a winning strategy or not, we must first check to see that its index is large enough. If so, we may use our rules from Lemma 6. If not, we will need to compute it directly, or look it up in a finite, pre-computed table.

We will now prove Theorem 2, by finishing the characterization for the remaining indices: $n\equiv 1, 3, \ldots 2d-1 \mod 4d$:

\begin{lemma}If $n\equiv 1, 3, 5, \ldots 2d-1 \mod 4d$, the base $2d$ representation of $n$ indicates whether or not the game $\{1, 2d, n\}$ is a first player win or loss. \end{lemma}

\begin{proof} Consider the block of even digits immediately preceding the final (least significant) odd digit of $n$. In the analogous Lemma 2, we said that if that block is even then $n$ is a loss,  and if it is odd then $n$ is a win. This was based on the idea that if we remove that block, we will have a smaller value of $n$, whose base $2d$ representation ended with 2 odd digits. We were guaranteed that this had no winning strategy. If the block we removed had even length, we would have the same situation, and if it was odd we would have a winning strategy (by handing the opponent the smaller value with no winning strategy).  

Again, the inductive structure of our argument remains, and only the base case has changed. We may still argue that if we remove the block of even numbers from the base $2d$ representation of $n$ we will get a smaller value $n^\prime$, and that the two states have the same relationship as before: if the block of even digits had an even length, they will have the same state (either both wins or both losses), whereas if the block removed had an odd length they will have the opposite state (if the smaller is a win the larger is a loss, and vice versa). 

What is different now is that, while we have a characterization for what happens when the base $2d$ representation ends in two odd digits, it's no longer the case that this is always a loss. To discover whether or not a particular $n$ is a win or loss we first remove the final block of even digits to obtain the smaller value $n^\prime$. By applying Lemma 6, we see that if $n^\prime$ is large enough this will be a loss, but we may have to look up whether this is a win or not in a finite lookup table. Either way, based on the (known) value at $n^\prime$ and the length of the block of even digits, we can easily compute the value at any $n$.
\end{proof}

As a remark, note that we cannot say that there is any threshold after which we can easily characterize the values of the indices where $n\equiv 1, 3, \ldots 2d-1 \mod 4d$. Consider the misere version of the simple game $\{1,2,n\}$. In this game, $n=3$ is a win, not a loss, but all other values when $n \equiv 3 \mod 4$ are losses. However, this single change causes an infinite number of other values with $n \equiv 1 \mod 4$ to change from wins to losses, or vice versa. All numbers that are one more than a power of 2, or equivalently, with binary representation $1,0 \ldots 0_k,1$ will depend on 3, and so be affected.

Lemmas 6 and 7 together give us a complete characterization of the games $\{a, 2d, n\}$, when $a'=1$, proving Theorem 2.

\section{Regularity of the Game Sequences}

A sequence $a_n$ is \emph{$k$-regular} if all residue classes modulo large powers of $k$ are formed by combining residue classes modulo smaller powers of $k$. More formally, we are interested in completely defining the sequence with recurrences of the following form:
$$a_{k^m+r} = \sum_{i=1}^L c_i a_{k^{m_i}+r_i},$$
where $m \geq m_i$, $0\leq r \leq k^m-1$, and $0\leq r_i \leq k^{m_i}-1$, and the $c_i$ can be any constants. Note that we do not require the smaller powers of $k$ to be equal to each other.

This property defines a large class of reasonably well-behaved and well understood sequences. A more in depth treatment can be found in \cite{Auto}.

We are particularly interested in the following related result: If a sequence is $k$-regular, and takes on only a finite number of values, it is $k$-automatic. This means that there is a discrete finite automaton that, when given the digits of $n$ in base $k$, will return the value of $a_n$. This has big implications for computability of these sequences. We will show that the Sprague-Grundy sequence of any game $\{1, 2d, n\}$ is $2d$-regular, by building up the family of recurrences that defines it. From this, we get that the sequence is also $2d$-automatic.

\subsection{Proof of Regularity}

\begin{theorem} The SG sequence of the game $\{1, 2d, n\}$ is $2d$-regular. \end{theorem}

\begin{proof} Let $n=R+2dc_1+4d^2c_2+8d^3c_3+16d^4c_4+\ldots$, for the moment with as many extra terms as we need. All the coefficients $c_i$ are integers in $[0,2d)$. This is equivalent to looking at the digits in the base-$2d$ representation of $n$, but we will write it as a polynomial so that if necessary we can borrow between terms. Note that because we have alternation, we only need to consider $R=0$ and $R=1$. 

We're looking for equivalences of the form $SG(n)=SG(r)$, for $r<n$. Of course, there are only 3 $SG$ values so there are many equations we could write of this form. But we would like to find large classes of these equivalences, where by only looking at the first few terms in the polynomial of $n$ (that is, the least significant digits in base $2d$) we can write down a formula for $r$ in terms of $n$. These classes of equalities will lead us to recurrences of the form $SG( (2d)^k)+i)=SG((2d)^m+j)$. In particular, we're interested in times when to write $r$ we have to delete digits from $n$, as this will give us equalities between different powers of $2d$.

To find these equivalences, we will use Lemma 1 and the following basic facts about mex: 
\begin{itemize}
\item $\text{mex}(0, 1)=2$, and $\text{mex}(0, 2)=1$.
\item (Corollary to above): $\text{mex}(0, \text{mex}(0,1))=1$, and $\text{mex}(0, \text{mex}(0,2))=2$.
\end{itemize}

For example, we can show that $SG_{1,2}(0+16c_4+\ldots)= SG_{1,2}(0+4c_4+\ldots)$:
\begin{align*}
SG_{1,2}(&0+ 16c_4+32c_5+\ldots)=\\
	&= \text{mex}\left(SG(-1+16c_4+32c_5+\ldots), SG(0+8c_4+16c_5+\ldots) \right)\\
	&= \text{mex}\left(0, SG(0+8c_4+16c_5+\ldots) \right) \qquad \text{Lemma 1}\\
	&= \text{mex}\left(0, \text{mex}\left(SG(-1+8c_4+16c_5+\ldots), SG(0+4c_4+8c_5+\ldots) \right) \right)\\
	&= \text{mex}\left(0, \text{mex}\left(0, SG(0+4c_4+8c_5+\ldots) \right) \right) \qquad \text{Lemma 1}\\
	&=SG(0+4c_4+8c_5+\ldots) \qquad \text{Third principle, above}
\end{align*}
Since we haven't used any terms above $c_4$, this tells us that the subsequence $SG_{1,2}(16n)$ equals the subsequence $SG_{1,2}(4n)$, term for term. It also tells us that if we add any other conditions on $c_4$ or coefficients of higher terms, we will still have equality. So, for example, it follows that $SG_{1,2}(32n)=SG_{1,2}(8n)$, $SG_{1,2}(32n+16)=SG_{1,2}(8n+4)$, etc. 

To this end, we will use a result for $n$ modulo $4d$ that we proved before:  If $n \equiv 2d+1, 2d+3, \ldots, 4d-1 \mod 4d$, or equivalently if $R=1$ and $c_1$ odd, then $SG(n)=0$. Based on this, we will break up our discussion of $n$ into the following 4 cases:  $R=0$ and $c_1$ even, $R=1$ and $c_1$ even, $R=0$ and $c_1$ odd, and $R=0$ and $c_1$ odd.

\subsubsection{R=1}

\noindent \textbf{Case 1:} $R=1$ and $c_1$ odd

In this case $SG(n)=0$, by Lemma (above).

\bigskip

\noindent \textbf{Case 2:} $R=1$,  $c_1$ even
$$SG(n) = \text{mex}\left( \vphantom{d^2} SG(2dc_1 +4d^2c_1 +\ldots), SG(1 + c_1 +2dc_2+ 4d^2c_3 + \ldots) \right)$$
\noindent If $c_2$ odd:
\begin{alignat*}{2}
SG&(1+2kc_1+ 4d^2c_2+8d^3c_3+\ldots) && \\
&=\text{mex}\left( SG(2dc_1 +4d^2c_2 +\ldots), SG(1 + c_1 +2dc_2 + 4d^2c_3 + \ldots)  \right)&&\\
 	&=\text{mex}\left( SG(2dc_1 +4d^2c_1 +\ldots), 0 \right) && 
	\\
	&= \text{mex}\left( \text{mex}\left( SG(-1 + 2dc_1 +4d^2c_2 +\ldots ), SG(c_1 +2dc_2 + \ldots)  \right), 0  \right)&&\\
	&= \text{mex}\left( \text{mex}\left( 0, SG(c_1 +2dc_2 \vphantom{d^2}+ \ldots)  \right), 0 \right). &&\end{alignat*}	
Since $SG(c_1 +2dc_2 + \ldots)$ is never zero, we have that $$SG(1+2dc_1+ 4d^2c_2+8d^3c_3+\ldots) = SG(c_1 +2dc_2 + \ldots).$$

\bigskip

\noindent If $c_2$ even,  $c_3$ odd:
\begin{align*}
SG(1+&2dc_1+4d^2c_2+8d^3c_3+\ldots)\\
	& = \text{mex} \left( SG(2dc_1 +4d^2c_2 +\ldots), SG(1 + c_1 +2dc_2 + 4d^2c_3\ldots) \right)\\
 	&=\text{mex} \left( \text{mex}\left( SG(-1 + 2dc_1 +4d^2c_2 + \ldots), SG(c_1 +2dc_2 + \ldots) \right) \right. ,\\
	&\qquad \qquad \quad \left. \text{mex}\left( SG(c_1 +2dc_2+ \ldots), SG(1+2c_2 + 4dc_3 + \ldots)\vphantom{d^2}\right) \right) \\
 	&=\text{mex} \left( \text{mex}\left( \vphantom{d^2}0, SG(2c_1 +4dc_2 + \ldots) \right), \text{mex}\left( \vphantom{d^2}SG(2c_1 +4dc_2+ \ldots), 0 \right) \right)\\
	&= 0
\end{align*}
Since $\text{mex}\left( \vphantom{d^2}0, SG(2c_1 +4dc_2 + \ldots) \right)\neq0$ and $\text{mex}\left( \vphantom{d^2}SG(2c_1 +4dc_2+ \ldots), 0 \right) \neq 0$, we have that the original is the mex of two non-zero elements, and so must always be zero.

If $c_2, c_3$ even, $c_4$ odd:
\begin{align*}
SG(1+ & 2dc_1+  4d^2c_2+8d^3c_3+16d^4c_4 + \ldots) \\
	&= \text{mex}\left( SG(2dc_1 +4d^2c_2 +\ldots), SG(1 + c_1 +2dc_2 +\ldots) \right)\\				
 	&=\text{mex}\left( \text{mex}\left( SG(-1 + 2dc_1 +4d^2c_2 + \ldots), SG(c_1+2dc_2 + \ldots) \right), \right.\\
	&\hspace{1 in} \left. SG(1 +2dc_2 + 4d^2c_3+ 8d^3c_4\ldots) \right)\\
	&=\text{mex} \left( \text{mex}\left(0, SG(c_1 +2dc_2 + \ldots) \vphantom{d^2} \right), SG(1 +2dc_2 + 4d^2c_3+ 8d^3c_4\ldots)  \right).  	  
	\end{align*}
But by the case above, we know that $ SG(1 +2dc_2 + 4d^2c_3+ 8d^3c_4\ldots)=0$, so 
\begin{align*}
SG(1+ & 2dc_1+  4d^2c_2+8d^3c_3+16d^4c_4 + \ldots) \\
	&=\text{mex} \left( \text{mex}\left(0, SG(c_1 +2dc_2 + \ldots) \vphantom{d^2} \right), SG(1 +2dc_2 + 4d^2c_3+ 8d^3c_4\ldots)  \right)  \\
	&=\text{mex} \left( \text{mex}\left(0, SG(c_1 +2dc_2 + \ldots) \vphantom{d^2} \right), 0 \right) . 	  	  
\end{align*}	
Since $SG(c_1 +2dc_2 + \ldots)$ is never zero (by assumption, $c_1$ even), we have that $$SG(1+2dc_1+ 4d^2c_2+8d^3c_3+\ldots) = SG(c_1 +2dc_2 + \ldots).$$

If $c_2, c_3, c_4$ even, we will break this down into two subcases: when $c_1=0$ and when $c_1\neq 0$. If $c_1=0$, $c_2, c_3, c_4$ all even:
\begin{align*}
SG&(1+4d^2c_2+8d^3c_3+16d^4c_4 + \ldots)\\
	&=\text{mex}\left( SG(4d^2c_2+8d^3c_3+16d^4c_4 + \ldots), SG(1+2dc_2 +4d^2c_3+8d^3c_4 + \ldots) \right)\\
	&=\text{mex}\left( \text{mex} \left( SG(-1+4d^2c_2+8d^3c_3+ \ldots), SG(2dc_2+ 4d^2c_3+8d^3c_4+ \ldots)\right), \right.\\
	&\hspace{0.5 in} \left. \text{mex}  \left( SG(2dc_2 +4d^2c_3+8d^3c_4 + \ldots), SG(1+2dc_3+4d^2c_4 +\ldots) \right) \right)\\
	&=\text{mex}\left( \text{mex} \left( 0, SG(2dc_2+ 4d^2c_3+8d^3c_4+ \ldots)\right), \right.\\
	&\hspace{0.5 in} \left. \text{mex}  \left( SG(2dc_2 +4d^2c_3+8d^3c_4 + \ldots), SG(1+2dc_3+4d^2c_4 +\ldots) \right) \right)
\end{align*}

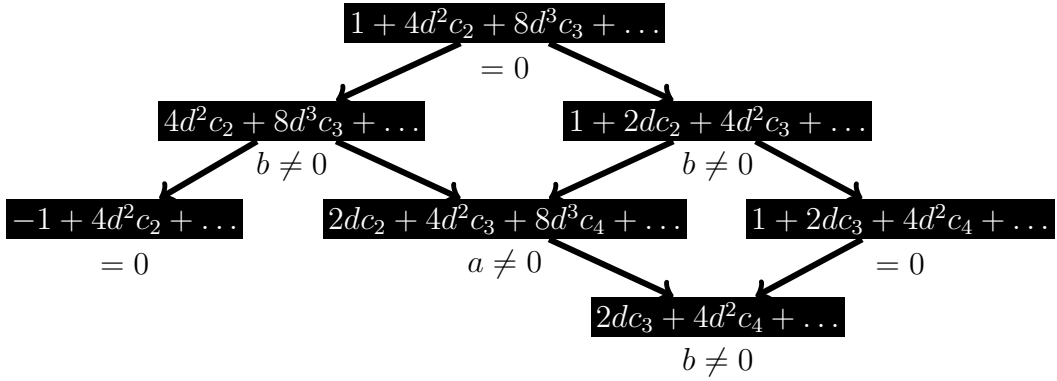
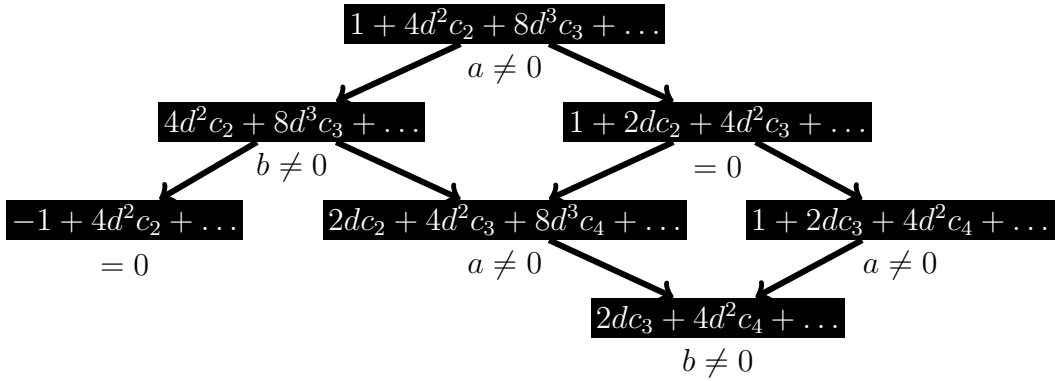
\begin{tikzpicture}
\tikzstyle{vert}=[rectangle,fill=black, text=white, minimum size =15, inner sep=0]
\tikzstyle{label}=[rectangle,fill=white, text=black, minimum size =15, inner sep=0]

\draw (0,0) node[vert] (Top){\,$1+4d^2c_2+8d^3c_3+ \ldots$\,};
\draw (-2.8,-1.3) node[vert] (L1) {\,$4d^2c_2+8d^3c_3 + \ldots$\,};
\draw (2.8,-1.3) node[vert] (R1){\,$1+2dc_2 +4d^2c_3+ \ldots$\,};
\draw (0,-2.6) node[vert] (RL){\,$2dc_2+ 4d^2c_3+8d^3c_4+ \ldots$\,};
\draw (5.2,-2.6) node[vert] (RR){\,$1+2dc_3+4d^2c_4 +\ldots$\,};
\draw (-5,-2.6) node[vert] (LL){\,$-1+4d^2c_2+ \ldots$\,};
\draw (2.8,-3.9) node[vert] (RRL){\,$2dc_3+4d^2c_4+\ldots$\,};

\draw (0,-0.6) node[label] (){$=0$};
\draw (2.8,-1.9) node[label] (){$b\neq 0$};
\draw (-2.8,-1.9) node[label] (){$b\neq 0$};
\draw (5.2,-3.2) node[label] (){$= 0$};
\draw (0,-3.2) node[label] (){$a \neq 0$};
\draw (-5,-3.2) node[label] (){$= 0$};
\draw (2.8,-4.5) node[label] (){$b\neq 0$};

\draw[->, line width=0.08cm] (Top) to (L1);
\draw[->, line width=0.08cm] (Top) to (R1);
\draw[->, line width=0.08cm] (R1) to (RL);
\draw[->,line width=0.08cm] (R1) to (RR);
\draw[->,line width=0.08cm] (L1) to (RL);
\draw[->,line width=0.08cm] (L1) to (LL);
\draw[->,line width=0.08cm] (RR) to (RRL);
\draw[->,line width=0.08cm] (RL) to (RRL);

\tikzstyle{vert}=[rectangle,fill=black, text=white, minimum size =15, inner sep=0, yshift=-5.5cm]
\tikzstyle{label}=[rectangle,fill=white, text=black, minimum size =15, inner sep=0, yshift=-5.5cm]

\draw (0,0) node[vert] (Top){\,$1+4d^2c_2+8d^3c_3+ \ldots$\,};
\draw (-2.8,-1.3) node[vert] (L1) {\,$4d^2c_2+8d^3c_3 + \ldots$\,};
\draw (2.8,-1.3) node[vert] (R1){\,$1+2dc_2 +4d^2c_3+ \ldots$\,};
\draw (0,-2.6) node[vert] (RL){\,$2dc_2+ 4d^2c_3+8d^3c_4+ \ldots$\,};
\draw (5.2,-2.6) node[vert] (RR){\,$1+2dc_3+4d^2c_4+\ldots $\,};
\draw (-5,-2.6) node[vert] (LL){\,$-1+4d^2c_2+ \ldots$\,};
\draw (2.8,-3.9) node[vert] (RRL){\,$2dc_3+4d^2c_4+\ldots$\,};

\draw (0,-0.6) node[label] (){$a \neq 0$};
\draw (2.8,-1.9) node[label] (){$=0$};
\draw (-2.8,-1.9) node[label] (){$b\neq0$};
\draw (5.2,-3.2) node[label] (){$a \neq 0$};
\draw (0,-3.2) node[label] (){$a \neq 0$};
\draw (-5,-3.2) node[label] (){$=0$};
\draw (2.8,-4.5) node[label] (){$b\neq 0$};

\draw[->, line width=0.08cm] (Top) to (L1);
\draw[->, line width=0.08cm] (Top) to (R1);
\draw[->, line width=0.08cm] (R1) to (RL);
\draw[->,line width=0.08cm] (R1) to (RR);
\draw[->,line width=0.08cm] (L1) to (RL);
\draw[->,line width=0.08cm] (L1) to (LL);
\draw[->,line width=0.08cm] (RR) to (RRL);
\draw[->,line width=0.08cm] (RL) to (RRL);

\draw (0,-11) node[rectangle, fill=white, text=black, text width=11cm](){When $c_1=0$, $c_2, \ldots c_4$ even, there are two possibilities for $SG(1+2d^3_4d^2c_4+\ldots)$ (found in the lower right). All the other $SG$ values can be determined from that value.};
\end{tikzpicture}

Notice that the term $SG(2dc_2 +4d^2c_3+8d^3c_4 + \ldots)$ shows up in both terms of the outer mex, and can never be zero. Further notice that  $SG(2dc_2 +4d^2c_3+8d^3c_4 + \ldots)$ and $SG(1+2dc_3+4d^2c_4 +\ldots)$ both depend on $SG(2dc_3+4d^2c_4+\ldots)$, which itself can never be zero. This means that when $SG(1+2dc_3+4d^2c_4 +\ldots)$ is non-zero, it must equal  $SG(2dc_2 +4d^2c_3+8d^3c_4 + \ldots)$. This is not true when $c_1\neq 0$.

If $SG(1+2dc_3+4d^2c_4 +\ldots)=0$, then both terms of the outer mex are given by $\text{mex}\left(0, SG(2dc_2 +4d^2c_3+8d^3c_4 + \ldots)\right)$, and are non-zero. So the original value must be zero as well. 

If $SG(1+2dc_3+4d^2c_4 +\ldots )=a \neq 0$, then it must be that $SG(2dc_2 +4d^2c_3+8d^3c_4 + \ldots)=a$, and so 
$SG(1+4d^2c_2+\ldots)=\text{mex}\left( \vphantom{d^2} \text{mex}(0,a), \text{mex}(a,a) \right) = \text{mex}\left( \vphantom{d^2} \text{mex}(0,a),0 \right)=a.$

In both cases, the value of the original is the same as $SG(1+2dc_3+4d^2c_4 +\ldots )$, so we have
$$SG(1+2dc_1+4d^2c_2+8d^3c_3+16d^4c_4 + \ldots)= SG(1 + 2dc_3 + 4d^2c_4 + \ldots).$$

If $c_1 \neq 0$, $c_2, c_3, c_4$ all even:
\begin{align*}
SG(& 1 +  2dc_1+4d^2c_2+8d^3c_3+16d^4c_4+ \ldots) \\
	&= \text{mex} \left( SG(2dc_1 +4d^2c_2 + 8d^3c_3 + \ldots), SG(1 +c_1 +2dc_2 + 4d^2c_3 +\ldots) \right)\\
 	&=\text{mex}\left( \text{mex}\left( SG(-1 + 2dc_1 +4d^2c_2 + \ldots), SG(c_1 +2dc_2 + 4d^2c_3 + \ldots) \right), \right. \\
	& \hspace{1in} \left.  \vphantom{d^2} \text{mex} \left( SG(c_1 +2dc_2 + 4d^2c_3 + \ldots), SG(1 + c_2 + 2dc_3 + \ldots) \right) \right) \\
 	&=\text{mex} \left( \text{mex} \left(0 , SG(c_1 +2dc_2 +4d^2c_3 + \ldots) \right), \right.\\
	& \hspace{1 in} \left.\text{mex}\left( SG(c_1+2dc_2+ 4d^2c_3 + \ldots), SG(1  +2dc_3 + \ldots) \right) \right)
\end{align*}
Note that the term $SG(c_1+2dc_2+ 4d^2c_3 + \ldots)$ can never be zero, and shows up in both terms of the outer mex. We also see that $SG(c_1+2dc_2+ 4d^2c_3 + \ldots)$ depends on the other term in the mex, $SG(1  +2dc_3 + \ldots)$ (here we are using the fact that $c_1\neq0$). From this we get that if $SG(1  +2dc_3 + \ldots)$ non-zero, it cannot have the same value as $SG(c_1+2dc_2+ 4d^2c_3 + \ldots)$: when one of them is 1, the other will be 2, and vice versa. 

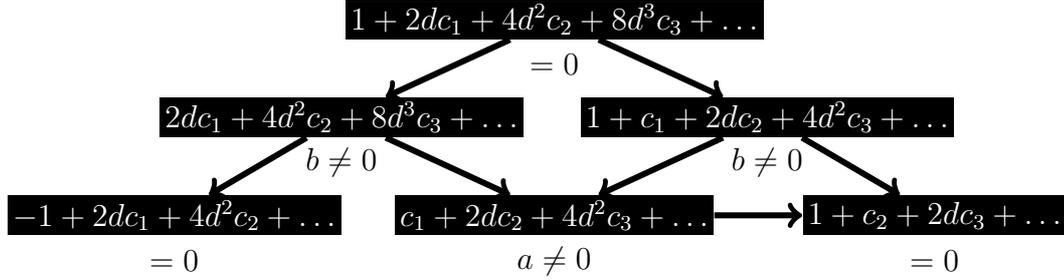
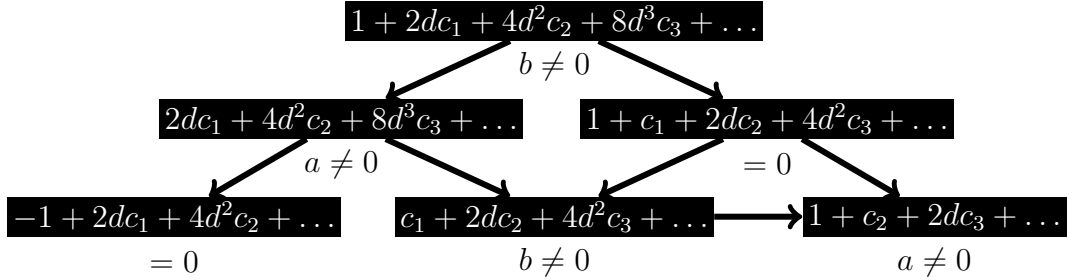
\begin{tikzpicture}
\tikzstyle{vert}=[rectangle,fill=black, text=white, minimum size =15, inner sep=0]
\tikzstyle{label}=[rectangle,fill=white, text=black, minimum size =15, inner sep=0]

\draw (0,0) node[vert] (Top){\,$1+2dc_1+4d^2c_2+8d^3c_3+ \ldots$\,};
\draw (-2.8,-1.3) node[vert] (L1) {\,$2dc_1+4d^2c_2+8d^3c_3 + \ldots$\,};
\draw (2.8,-1.3) node[vert] (R1){\,$1+c_1+2dc_2 +4d^2c_3+ \ldots$\,};
\draw (0,-2.6) node[vert] (RL){\,$c_1+2dc_2+ 4d^2c_3+ \ldots$\,};
\draw (5,-2.6) node[vert] (RR){\,$1+c_2+2dc_3 +\ldots$\,};
\draw (-5,-2.6) node[vert] (LL){\,$-1+2dc_1+4d^2c_2+ \ldots$\,};

\draw (0,-0.6) node[label] (){$=0$};
\draw (2.8,-1.9) node[label] (){$b\neq 0$};
\draw (-2.8,-1.9) node[label] (){$b\neq 0$};
\draw (5,-3.2) node[label] (){$= 0$};
\draw (0,-3.2) node[label] (){$a \neq 0$};
\draw (-5,-3.2) node[label] (){$= 0$};

\draw[->, line width=0.08cm] (Top) to (L1);
\draw[->, line width=0.08cm] (Top) to (R1);
\draw[->, line width=0.08cm] (R1) to (RL);
\draw[->,line width=0.08cm] (R1) to (RR);
\draw[->,line width=0.08cm] (L1) to (RL);
\draw[->,line width=0.08cm] (L1) to (LL);
\draw[->,line width=0.08cm] (RL) to (RR);


\tikzstyle{vert}=[rectangle,fill=black, text=white, minimum size =15, inner sep=0, yshift=-4.5cm]
\tikzstyle{label}=[rectangle,fill=white, text=black, minimum size =15, inner sep=0, yshift=-4.5cm]

\draw (0,0) node[vert] (Top){\,$1+2dc_1+4d^2c_2+8d^3c_3+ \ldots$\,};
\draw (-2.8,-1.3) node[vert] (L1) {\,$2dc_1+4d^2c_2+8d^3c_3 + \ldots$\,};
\draw (2.8,-1.3) node[vert] (R1){\,$1+c_1+2dc_2 +4d^2c_3+ \ldots$\,};
\draw (0,-2.6) node[vert] (RL){\,$c_1+2dc_2+ 4d^2c_3+ \ldots$\,};
\draw (5,-2.6) node[vert] (RR){\,$1+c_2+2dc_3 +\ldots$\,};
\draw (-5,-2.6) node[vert] (LL){\,$-1+2dc_1+4d^2c_2+ \ldots$\,};

\draw (0,-0.6) node[label] (){$b \neq 0$};
\draw (2.8,-1.9) node[label] (){$=0$};
\draw (-2.8,-1.9) node[label] (){$a \neq0$};
\draw (5,-3.2) node[label] (){$a \neq 0$};
\draw (0,-3.2) node[label] (){$b \neq 0$};
\draw (-5,-3.2) node[label] (){$=0$};

\draw[->, line width=0.08cm] (Top) to (L1);
\draw[->, line width=0.08cm] (Top) to (R1);
\draw[->, line width=0.08cm] (R1) to (RL);
\draw[->,line width=0.08cm] (R1) to (RR);
\draw[->,line width=0.08cm] (L1) to (RL);
\draw[->,line width=0.08cm] (L1) to (LL);
\draw[->,line width=0.08cm] (RL) to (RR);


\draw (0,-9) node[rectangle, fill=white, text=black, text width=10.7cm](){When $c_1\neq0$, $c_2,\ldots c_4$ even, there are two possibilities for $SG(1+2dc_3+4d^2c_4+\ldots)$ (lower right), and all other values can determined from that.};
\end{tikzpicture}

We can show that there is a nice relationship between $SG(1+2dc_1+4d^2c_2+8d^3c_3+16d^4c_4+ \ldots)$ and $SG(1 +c_2 +2dc_3 + 4d^2c_4 \ldots)$. When  $SG(1 +c_2 +2dc_3 + 4d^2c_4 \ldots)=0$, then $SG(1+2dc_1+4d^2c_2+8d^3c_3+16d^4c_4+ \ldots)$ is the mex over two non-zero elements, and so must be zero as well.  

But when $SG(1 +c_2 +2dc_3 + 4d^2c_4 \ldots)=a \neq 0$, we must have that $SG(c_1+2dc_2+ 4d^2c_3 + \ldots)= b\neq0$. In this case,
\begin{align*}
SG(& 1 +  2dc_1+4d^2c_2+8d^3c_3+16d^4c_4+ \ldots) \\
&=\text{mex} \left( \text{mex} \left(0 , SG(c_1 +2dc_2 +4d^2c_3 + \ldots) \right), \right.\\
	& \hspace{1 in} \left.\text{mex}\left( SG(c_1+2dc_2+ 4d^2c_3 + \ldots), SG(1  +2dc_3 + \ldots) \right) \right)\\
	&= \text{mex} \left( \text{mex} \left(0 , b \right), \text{mex}\left( b, a) \right) \right)\\
	&= \text{mex} \left( \text{mex} \left(0 , b \right), 0 \right) =b
\end{align*}

So we have that $SG(1 +  2dc_1+4d^2c_2+8d^3c_3+16d^4c_4+ \ldots)=SG(1 +c_2 +2dc_3 + 4d^2c_4 \ldots)$ when they are both zero, and they are opposites (one is 1 while the other is 2) when they are both not zero.

We create a new function $SG^*(n)$ by assigning $SG^*(n)=(3-SG(n)) \mod 3$. In this way, $SG^*(n)$ and  $SG(n)$ are either both zero, or both non-zero and not equal. Note that $SG^{**}(n)=SG(n)$. We call this new function the starred value of $SG(n)$, or the opposite value if we are looking at a class of values $n$ where $SG(n)\neq 0$. Using this new function, we can summarize what we have learned with the following: If $c_1\neq 0$, $c_2, c_3, c_4$ even, then
$$SG(1 +  2dc_1+4d^2c_2+8d^3c_3+16d^4c_4+ \ldots)=SG^*(1 +2dc_3 + 4d^2c_4 \ldots).$$
We have 5 reductions available for $SG(1 +  2dc_1+4d^2c_2+8d^3c_3+16d^4c_4+ \ldots)$, summarized in this table:

\begin{center}
\begin{tabular}{|lll|}
\hline
&Reduction of the original& Conditions\\
\hline
1)& $SG(c_1 +2dc_2 + 4d^2c_3+\ldots)$ & $c_2$ odd \\
2)& 0 & $c_2$ even, $c_3$ odd \\
3)&$SG(c_1 +2dc_2 + 4d^2c_3 +\ldots)$  & $c_2, c_3$ even, $c_4$ odd \\
4)&$SG(1 + 2dc_3 + 4d^2c_4 + \ldots)$ & $c_1= 0$, $c_2, c_3, c_4$ even\\
5)&$ SG^*(1 +2dc_3 + 4d^2c_4+ \ldots)$ & $c_1\neq 0$, $c_2, c_3, c_4$ even\\
\hline
\end{tabular}
\end{center}

In Rule 5, we don't get an equality with $SG$ of another value, but with it's starred counterpart. Notice, however, that the new value still satisfies the constraints of this class of digits, namely that $R=1$ and that $c_3$, coefficient of $2d$, is even. So we will use the rules in this table on our smaller number, $1 +2dc_3 + 4d^2c_4 + \ldots$, to find a further reduction.

In fact, Rule 1 does not apply, since $c_4$ even by assumption. If $c_5$ is odd, we apply Rule 2 which tells us that 
\begin{align*}
SG(1 +  2dc_1+4d^2c_2 &+8d^3c_3+16d^4c_4+ \ldots)\\
	&=SG^*(1 +2dc_3 + 4d^2c_4 +\ldots)\\
	& = 0.
\end{align*}
  
If $c_5$ even but $c_6$ odd, we are in the case of Rule 3 and when we apply it we get 
\begin{align*}
SG(1 +  2dc_1+4d^2c_2+&8d^3c_3+16d^4c_4+ \ldots)\\
&= SG^*(1 +2dc_3 + 4d^2c_4 +\ldots) \\
&= SG^*(c_3 +2dc_4 +4d^2c_5 + \ldots).
\end{align*}
$SG(c_3 +2dc_4 +4d^2c_5 + \ldots)$ is never zero (Lemma 1 applies, since $c_3$ even, by assumption), so $SG^*(c_3 +2dc_4 +4d^2c_5 + \ldots)=\text{mex}\left(0, SG(c_3 +2dc_4 +4d^2c_5 + \ldots) \right)$. 
From this we see that $SG(2dc_3+4d^2c_4+ \ldots)$ also equals $SG^*(c_3 +2dc_4 +4d^2c_5 + \ldots)$. Therefore, it must be that $$SG(1 +  2dc_1+4d^2c_2+8d^3c_3+16d^4c_4+ \ldots) = SG(2dc_3+4d^2c_4+ \ldots).$$

If $c_3=0$, $c_4, c_5, c_6$ all even, then we are in the case of Rule 4.  
\begin{align*}
SG(1 +  2dc_1+4d^2c_2+&8d^3c_3+16d^4c_4+ \ldots)\\
&= SG^*(1 +2dc_3 + 4d^2c_4 +\ldots)\\
&= SG^*(1+ 2dc_5+ 4d^2c_6 + \ldots)
\end{align*} 
We are looking for some other value that causes us to reduce in such a way that we get equality with a starred value, and the values for $c_1$ and $c_2$ do just that.
Consider taking away the $c_3$ and $c_4$, so that we are left with $SG(1 +  2dc_1 +4d^2c_2+8d^3c_5+16d^4c_6+ \ldots)$. Rule 5 applies (since by assumption $c_1\neq0$, $c_2, c_5, c_6$ all even), and so it must equal $SG^*(1+ 2dc_5+ 4d^2c_6 + \ldots)$. This gives us that
$$SG(1 +  2dc_1+4d^2c_2+8d^3c_3+16d^4c_4+ \ldots)= SG(1 +  2dc_1 +4d^2c_2+8d^3c_5+16d^4c_6+ \ldots).$$

Finally, in the case of Rule 5, we have
\begin{align*}
SG(1 +  2dc_1+4d^2c_2+&8d^3c_3+16d^4c_4+ \ldots)\\
 &= SG^*(1 +2dc_3 + 4d^2c_4 + \ldots) \\
 &= SG^{**}(1+ 2dc_5+ 4d^2c_6 + \ldots)\\ 
 &= SG(1+ 2dc_5+ 4d^2c_6 + \ldots).
\end{align*}

To aid with further reductions in the other cases, we summarize all our reductions here. The final column tells us the shift, or how many orders of $2d$ we lose in our reduction.
\begin{center}
\begin{tabular}{|l|}
\hline
Reductions for $SG(1 + 2dc_1 +4d^2c_2+8d^3c_3+ \ldots)$ when $c_1$ even:\\
\hline
\begin{tabular}{lllc}
1)& $SG(c_1 +2dc_2 + 4d^2c_3+\ldots)$ & $c_2$ odd & 1 \\
2)& 0 & $c_2$ even, $c_3$ odd & NA \\
3)&$SG(c_1 +2dc_2 + 4d^2c_3 +\ldots)$  & $c_2, c_3$ even, $c_4$ odd & 1 \\
4)&$SG(1 + 2dc_3 + 4d^2c_4 + \ldots)$ & $c_1= 0$, $c_2, c_3, c_4$ even & 2\\
5.1)& 0 & $c_1\neq 0$, $c_2, c_3, c_4$ even, $c_5$ odd & NA\\
5.2)&$ SG(2dc_3 + 4d^2c_4+ \ldots)$ & $c_1\neq 0$, $c_2, c_3, c_4, c_5$ even, $c_6$ odd & 2\\
5.3)&$ SG(1 +2dc_1 + 4d^2c_2+8d^3c_5+ \ldots)$ & $c_1\neq 0$, $c_2, c_4, c_5, c_6$ even, $c_3=0$ &2\\
5.4)&$ SG(1 +2dc_5 + 4d^2c_6+ \ldots)$ & $c_1, c_3 \neq 0$, $c_2, c_3, c_4, c_5, c_6$ even& 4\\
\end{tabular}\\
\hline
\end{tabular}
\end{center}

\bigskip

\noindent \textbf{Case 3:} $R=0$, $c_1$ odd
\begin{align*}
SG&(2dc_1+4d^2c_2+8d^3c_3+16d^4c_4+\ldots) \\
&=\text{mex}\left( SG(1+2d(c_1-1)+4d^2c_2+8d^3c_3+\ldots), SG(1+2dc_2+4d^2c_3+\ldots) \right) 
\end{align*}

Note that instead of subtracting 1 in the standard way, we have borrowed from $c_1$. This is possible because $c_1$ odd, therefore $c_1 \geq 1$. Note also that this now puts the first term in the framework of Case 2, above. We will denote the first term in the mex by $LHS$, and the second term by $RHS$.

In many cases, we can apply some of the rules for the previous case to get one of the sides to disappear entirely:
\begin{itemize}
\item If $c_2$ odd, then $RHS=0$, by Lemma. Reduce $LHS$ with Rule 1.

\item If $c_2$ even, $c_3$ odd, then $LHS=0$, by Rule 2. Reduce $RHS$ with Rule 1.

\item If $c_2, c_3$ even, $c_4$ odd, then $RHS=0$, by Rule 2. Reduce $LHS$ with Rule 3.

\item If $c_1 \neq 1$, $c_2, c_3, c_4$ even, $c_5$ odd, then $LHS=0$, by Rule 5.1.  Reduce $RHS$ with Rule 3.

\item If $c_1=1$, then $LHS=SG(1+2dc_3+4d^2c_4+8d^3c_5+\ldots)=0$ by Rule 2. Reduce $RHS$ with Rule 3.

\item If $c_2 \neq 0$, $c_2, c_3, c_4, c_5$ even, $c_6$ odd, then $RHS=0$, by Rule 5.1. If $c_1\neq 1$, reduce $LHS$ by Rule 5.2. If $c_1=1$, reduce $LHS$ by Rule 4, and then Rule 3. 

\item If $c_2=0$, then $RHS=SG(1+2dc_4+4d^2c_5+8d^3c_6+\ldots)=0$ by Rule 2. If $c_1\neq 1$, reduce $LHS$ by Rule 5.2. If $c_1=1$, reduce $LHS$ by Rule 4, and then Rule 3.

\end{itemize}

In each case above, we have that $SG(2dc_1+4d^2c_2+8d^3c_3+16d^4c_4+\ldots)$ is the mex of another $SG$ value, and of 0. Our plan, for each of these cases, is to use the table above to reduce the remaining non-zero side, and then find some other number that depends on 0 and the reduced value. 

The final rule used to reduce the remaining side is always Rule 1, 3 or 5.1. These have the feature that the resulting reduced value is always even. When the resulting value is even, finding another value that depends on 0 and the reduced value is not difficult. If $r$ is the reduced value, then $SG(2d \cdot (r)) = \text{mex} (0, r)$. The only thing to check is that even if we multiply by $2d$ we still have a reduction. 

If the net effect of the reduction is a change in order of at least two for the LHS, or at least one for the RHS, then even if we multiply by $2d$ the result will be a value of a smaller order which we know must be a reduction. So the only times to be concerned are when we reduce the LHS with Rule 1 or Rule 3, which only involve a shift of one order. 

When we do the reduction of $LHS=SG(1+2d(c_1-1)+4d^2c_2+8d^3c_3+\ldots)$ with Rule 1 (or Rule 3, equivalently), we get $SG(c_1-1+2d(c_2)+4d^2c_3+8d^3c_4+\ldots)$. Now when we multiply through by $2d$, we get $SG(2d(c_1-1)+4d^2c_2+8d^3c_3+\ldots)$. This is strictly less than what we started with: $SG(2dc_1+4d^2c_2+8d^3c_3+\ldots)$, since we have decreased the value $c_1$ by one. So this is truly a reduction, even though there is no loss of order.

\subsubsection{Summary of $R=1$, $c_1$ even}
Table summary of results from above:

\begin{center}
\begin{tabular}{|l|l|c|}
\hline
Reduction of orginal & Conditions & Shift\\
\hline
$SG(2d(c_1-1)+4d^2c_2+\ldots)$& $c_2$ odd &0\\
$SG(2dc_2+4d^2c_3+8d^3c_4+\ldots)$&$c_2$ even, $c_3$ odd &1\\
$SG(2d(c_1-1)+4d^2c_2+\ldots)$& $c_2, c_3$ even, $c_4$ odd& 0\\
$SG(2dc_2+4d^2c_3+8d^3c_4+\ldots)$& $c_2, c_3, c_4$ even, $c_5$ odd&1\\
$SG(4d^2c_3+8d^3c_4+\ldots)$&  $c_1 \neq 1$, $c_2, c_3, c_4, c_5$ even, $c_6$ odd&1\\
$SG(2dc_3+4d^2c_4+\ldots)$ & $c_1 = 1$, $c_2, c_3, c_4, c_5$ even, $c_6$ odd&2 \\
\hline
\end{tabular}
\end{center}

\subsubsection{R=0, $c_1$ odd}

The harder cases to handle are when $c_2 \ldots c_6$ all even. We don't get one of the two terms in the mex going to zero. Instead, each reduction we make will reduce the order by 2, and will either be an equality (Rule 4), or an opposite (Rule 5).  We will hope to find a reduction of  $LHS$ and a reduction of $RHS$ so that their orders are within one, and they are either both equal or both opposite. This will cover almost all the cases. Here is a schematic for when we can do this:

\begin{center}
\begin{tabular}{c|c|c|c|c}
& 
\begin{tabular}{c}
$c_1=1$\\
\end{tabular}&
\begin{tabular}{c}
$c_1\neq1$, $c_2 \neq 0$\\
\end{tabular}&
\begin{tabular}{c}
$c_1\neq 1$,
$c_2= 0$\\
$c_3 \neq 0$
\\
\end{tabular}&
\begin{tabular}{c}
$c_1\neq 1$,
$c_2= 0$\\
$c_3 = 0$,
$c_4\neq0$\\
\end{tabular}\\
\hline
Shift: & LHS \quad RHS & LHS \quad RHS & LHS \quad RHS & LHS \quad RHS \\
\begin{tabular}{r}
0\\
1\\
2\\
3\\
4\\
5\\
\end{tabular}&
\begin{tabular}{ll}
SG& \\
 & \textbf{SG}\\
\textbf{SG}& \\
&\\
&\\
&\\
\end{tabular}&
\begin{tabular}{ll}
SG& \\
 & SG\\
\textbf{SG}$^*$& \\
&\textbf{SG}$^*$\\
&\\
&\\
\end{tabular}&
\begin{tabular}{ll}
SG& \\
 & SG\\
SG$^*$& \\
&\textbf{SG}\\
\textbf{SG}&\\
&\\
\end{tabular}&
\begin{tabular}{ll}
SG& \\
 & SG\\
SG$^*$& \\
&SG\\
\textbf{SG}$^*$&\\
& \textbf{SG}$^*$\\
\end{tabular}\\
\end{tabular}
\end{center}

To do really argue this, we need to show that when they are reduced at least by one we can find a reduction.

\begin{lemma}
If there are reductions of at least one of LHS and RHS that bring them to within 1 order of each other, and the reductions are to equality, then there is a reduction of the original.
\end{lemma}

\begin{proof}
The way the reductions go, we will have one of the sides be $SG(1+2dc_i+4d^2c_{i+1}+\ldots)$, and the other side will be $SG(1+2dc_{i+1}+4d^2c_{i+2} + \ldots)$. But since at least one of the sides has been reduced, we know that $i\geq 2$. We also know that in the cases we're considering, $c_i\in\{c_2, \ldots c_5\}$, and therefore is even. Consider the value $SG(2d(c_i+1)+4d^2c_{i+1}+\ldots)$. By definition, this equals $\text{mex}\left(SG(1+2dc_i+4d^2c_{i+1}+\ldots),  SG(1+2dc_{i+1}+4d^2c_{i+2} + \ldots)   \right)$, and so must equal the original. However, because $i\geq 2$, this is of a lower order. Even though we have raised one of the coefficients, the lower order ensures that it's still a reduction.
\end{proof}

\begin{lemma}
If there are reductions of at least one of LHS and RHS that bring them to within 1 order of each other, and the reductions are to opposites, then there is a reduction of the original.
\end{lemma}

\begin{proof}
The reductions will be of the form $SG^*(1+2dc_i+4d^2c_{i+1}+\ldots)$, and $SG^*(1+2dc_{i+1}+4d^2c_{i+2} + \ldots)$. We know that since the original is never zero, the LHS and the RHS are never both non-zero. Moreover, the larger of the two depends on the smaller, so they can never both be zero. Thus, it follows that their support is exactly complimentary.

Consider $\text{mex}\left(SG^*(1+2dc_i+4d^2c_{i+1}+\ldots), SG^*(1+2dc_{i+1}+4d^2c_{i+2} + \ldots)\right)$, and how this relates to $SG(2d(c_i+1)+4d^2c_{i+1}+\ldots)$, which is the mex of the same two terms but without the $^*$ in each. One of the two terms in the mex will always be zero, the other will be opposite what we would get if we had equality instead of opposite. So our result will always be opposite what we would have expected, and overall we have
\begin{align*}
SG^*(&2d(c_i+1)+4d^2c_{i+1}+\ldots)\\
&=\text{mex}\left(SG^*(1+2dc_i+4d^2c_{i+1}+\ldots), SG^*(1+2dc_{i+1}+4d^2c_{i+2} + \ldots)\right)
\end{align*}

This means we have that our original $SG(2dc_1+4d^2c_2+8d^3c_3+\ldots)$ equals $SG^*(2d(c_i+1)+4d^2c_{i+1}+\ldots)$, but this is still not a reduction. To fix this, we will multiply by $2d$:
$$SG(2dc_1+4d^2c_2+8d^3c_3+\ldots)=SG(4d^2(c_i+1)+8d^3c_{i+1}+\ldots).$$
We check that, by definition 
\begin{align*}
SG&(4d^2(c_i+1)+8d^3c_{i+1}+\ldots)=\\
	&=\text{mex} \left(SG(-1+4d^2(c_i+1)+8d^3c_{i+1}+\ldots),SG(2d(c_i+1)+4d^2c_{i+1}+\ldots)\right)\\
	&=\text{mex} \left(0,SG(2d(c_i+1)+4d^2c_{i+1}+\ldots)\right) \qquad \text{ Lemma 1}
\end{align*}
Since $SG(2d(c_i+1)+4d^2c_{i+1}+\ldots)$ is never zero, if we take its mex with 0, we do get $SG(4d^2(c_i+1)+8d^3c_{i+1}+\ldots)=SG^*(2d(c_i+1)+4d^2c_{i+1}+\ldots)$.

Finally, we must check that this is, in fact, a reduction. To get reductions to opposites, rather than equality, we must have reduced both the LHS and the RHS at least once each. Therefore, $c_i\geq3$. From this we see that we have reduced our order, and so this is a reduction.
\end{proof}

The single remaining case to consider is $c_1\neq1, c_2= c_3 = c_4=0$:\\
\begin{center}
\begin{tabular}{r|ll}
Shift: & LHS & RHS\\
\hline
0&SG& \\
1& & SG\\
2&SG$^*$& \\
3&&SG\\
4&SG$^*$&\\
5&& SG\\
\end{tabular}
\end{center}
We cannot find a pair that are reduced, within one order, and both equal or both opposites, so Lemmas 6 and 7 don't help us. Instead, we will use the fact that the 5th and 6th rows of the chart are the same as the 3rd and 4th rows.

After applying the reductions, we have the following:
\begin{align*}
LHS &= SG(1+2d(c_1-1)+4d^2c_2+8d^3c_3+ 16d^4c_4+32d^5c_5+64d^6c_6+\ldots)\\
  	&= SG^*(1+2dc_3+4d^2c_4+8d^3c_5+ 16d^4c_6+\ldots) \qquad\text{ Rule 5}\\
	&= SG^*(1+2dc_5+4d^2c_6+\ldots) \qquad \text{Rule 4}
\end{align*}

\begin{align*}
RHS &= SG(1+2dc_2+4d^2c_3+8d^3c_4+ 16d^4c_5+32d^5c_6+\ldots)\\
  	&= SG(1+2dc_4+4d^2c_5+8d^3c_6+\ldots) \qquad\text{ Rule 4}\\
	&= SG(1+2dc_6+\ldots) \qquad \text{Rule 4}
\end{align*}
We would like to find a smaller value, in the same class, that when we apply the definition and the simplification rules, we get the same awkward split into a term on the left that is opposite, and a term on the right that is equal. The first few values of the original value gave us that exact split, so we will use them again:
\begin{align*}
SG&(2dc_1+4d^2c_2+8d^3c_5+16d^4c_6+\ldots)\\
	&=\text{mex}\left(SG(1+2d(c_1-1) +4d^2c_2+8d^3c_5 + 16d^4c_6 ), SG(1+2dc_2+4d^2c_5+ 8d^3c_6+ \ldots) \right)\\
	&=\text{mex}\left(SG^*(1+2dc_5 +4d^2c_6+\dots ), SG(1+2dc_6+ \ldots) \right) \qquad \text{ Rules 4 and 5}
\end{align*}
So we have that $$SG(2dc_1+4d^2c_2+8d^3c_3+16d^4c_4+\ldots)   = SG(2dc_1+4d^2c_2+8d^3c_5+16d^4c_6+\ldots).$$

Tabular summary of all the reductions, for $R=0$, $c_1$ odd:
\begin{center}
\begin{tabular}{|l|l|c|}
\hline
Reduction of orginal & Conditions & Shift\\
\hline
$SG(2d(c_1-1)+4d^2c_2+\ldots)$& $c_2$ odd &0\\
$SG(2dc_2+4d^2c_3+8d^3c_4+\ldots)$&$c_2$ even, $c_3$ odd &1\\
$SG(2d(c_1-1)+4d^2c_2+\ldots)$& $c_2, c_3$ even, $c_4$ odd& 0\\
$SG(2dc_2+4d^2c_3+8d^3c_4+\ldots)$& $c_2, c_3, c_4$ even, $c_5$ odd&1\\
$SG(4d^2c_3+8d^3c_4+\ldots)$&  $c_1 \neq 1$, $c_2, c_3, c_4, c_5$ even, $c_6$ odd&1\\
$SG(2dc_3+4d^2c_4+\ldots)$ & $c_1 = 1$, $c_2, c_3, c_4, c_5$ even, $c_6$ odd&2 \\
\hline
& $c_2, \ldots c_6$ all even, and:&\\
$SG(2d(c_2+1)+4d^2c_3+\ldots)$& $c_1=1$ &1\\
$SG(4d^2(c_3+1)+8d^3c_4+\ldots)$& $c_1\neq 1$, $c_2\neq0$ &1\\
$SG(2d(c_4+1)+4d^2c_5+\ldots)$&  $c_1\neq 1$, $c_2=0$, $c_3\neq0$&3\\
$SG(4d^2(c_5+1)+8d^3c_6+\ldots)$&  $c_1\neq 1$, $c_2=c_3=0$, $c_4\neq0$&3\\
$SG(2dc_1+8d^3c_5+16d^4c_6+\ldots)$& $c_1\neq 1$, $c_2=c_3=c_4=0$&2\\
\hline
\end{tabular}
\end{center}

\bigskip

\subsubsection{$R=0$, $c_1$ even}

\noindent \textbf{Case 4:} $R=0$, $c_1$ even

\bigskip

Having done so much work tabulating the results in Cases 2 and 3, Case 4 has pleasantly few cases.
\begin{align*}
SG&(2dc_1+4d^2c_2+8d^3c_3+16d^4c_4+\ldots\\
	&=\text{mex}\left( SG(-1+2dc_1+4d^2c_2+\ldots), SG(c_1+2dc_2+4d^2c_3+\ldots)\right)
\end{align*}
Assuming that there is some non-zero coefficient in the LHS (i.e. that the value is at least as big as $2d$), we can borrow from one of the terms and we will get by Lemma 1 that the LHS is always zero. Therefore, we have
\begin{align*}
SG(2dc_1+4d^2c_2+8d^3c_3+16d^4c_4+\ldots)
	=\text{mex}\left(0, SG(c_1+2dc_2+4d^2c_3+\ldots)\right)
\end{align*}

Because $c_1$ even, the RHS must fall into either Case 3, or back into Case 4 (depending on the parity of $c_2$, and whether or not $c_1$=0). If the RHS falls back into Case 4, it must be that either $c_1=0$ and $c_2$ even, or that $c_1>0$ and $c_2$ is odd. In this second case, the principle of alternation tells us that RHS$=SG(2d(c_2+1) +4d^2c_3+ \ldots)$. But it's possible that $c_2=2d-1$, and we would have to carry over into larger coefficients. To be consistent, in either case we will rewrite it as $(SG(2dc_2^\prime+4d^2c_3^\prime+\ldots)$, where $c_2^\prime$ is even. 
\begin{align*}
SG&(2dc_1+4d^2c_2+8d^3c_3+16d^4c_4+\ldots)\\
	&=\text{mex}\left(0, SG(2dc_2^\prime+4d^2c_3^\prime+\ldots)\right)\\
	&=\text{mex}\left(0, \text{mex}\left(SG(-1+2dc_2^\prime+4d^2c_3^\prime+\ldots), SG(c_2^\prime+2dc_3^\prime+4d^2c_4^\prime+\ldots)\right)\right)\\
	&=\text{mex}\left(0, \text{mex}\left( 0, SG(c_2^\prime+2dc_3^\prime+4d^2c_4^\prime+\ldots)\right)\right)
\end{align*}
Since $SG(c_2^\prime+2dc_3^\prime+4d^2c_4^\prime+\ldots)$ is never zero (by Lemma 1), we have that $$SG(2dc_1+4d^2c_2+8d^3c_3+16d^4c_4+\ldots)=SG(c_2^\prime+2dc_3^\prime+4d^2c_4^\prime+\ldots).$$
This is a proper reduction, even though it's possible that the coefficients have increased, because the order has decreased.

If instead $c_1=0$ and $c_2$ odd, or $c_1>0$ and $c_2$ even, then the RHS falls into Case 3. In the second case, using the principles of alternation, we would write the RHS as $SG(2d(c_2+1)+ 4d^2c_3+\ldots)$, and there is no danger of carrying since $c_2<2d-1$. Either way, from here we reduce RHS using one of the rules of Case 3. Suppose the reduced value is $r$, then we will have that $SG(2dc_1+4d^2c_2+8d^3c_3+16d^4c_4+\ldots) = SG(2d\cdot r)$.

The final to thing to check is whether or not this is a proper reduction. In the first case, it must be because we haven't increased the RHS - it is exactly the original divided by $2d$. By the rules in Case 3, $r<RHS$, and so $2dr<2d\cdot RHS$. 

In the second case, this argument fails because we have actually artificially increased the RHS to get it into the right form. If $r$ is a reduction with shift zero, while it may be a proper reduction of the augmented RHS, it could potentially be greater than the original divided by $2d$. But in fact there is only one reduction with shift zero, and it takes $SG(2d(c_2+1)+ 4d^2c_3+\ldots)$ to $SG(2dc_2+ 4d^2c_3+\ldots)$. So even when we multiply again by $2d$, we will be missing the $c_1$ term, which we know is non-zero in this case. 
\end{proof}

\bibliographystyle{amsplain}

\begin{thebibliography}{99}
\frenchspacing
\bibitem{Fraenkel1}
A.S. Fraenkel, {\it The vile, dopey, evil and odious game players}, Discrete
Math. 312 (2012) 42-46, special volume in honor of the 80th birthday
of Gert Sabidussi.

\bibitem{Fraenkel2}
A.S. Fraenkel, {\it Aperiodic subtraction games}, Electronic J.
Combinatorics vol. 18(2) P19 12pp., 2011.

\bibitem{Guo}
A. Guo, {\it Winning strategies for aperiodic subtraction games,} ArXiV, \url{http://arxiv.org/abs/1108.1239}

\bibitem{Auto}
J. Allouche and J. Shallit. {\it Automatic Sequences: Theory, Applications, Generalizations}. Cambridge University Press, 2003.

\bibitem{first} E. Berlekamp and J. P. Buhler. \textit{Puzzles Column}, \textit{Emissary -- MSRI Gazette.} p. 6, Fall 2009. 


\end{thebibliography}
{

}

\end{document}